\documentclass[11pt]{article}

\usepackage{amssymb,amsmath,amsthm,color,graphics,graphicx} 

\begin{document}

\newcommand{\firstorsecond}{first}
\newcommand{\bloom}[1]{(\emph{Bloom's level {#1}})}
\newcommand{\listkeywords}[1]{\noindent \textbf{Keywords:} #1}
\newcommand{\assignment}[3]{\textbf{{\small{SAMPLE}} Assignment \#{#1}, given Week {#2}}}
\newcommand{\fordeadline}[1]{\noindent \emph{Due at class time\ {#1}\ from today:}}
\newcommand{\discuss}{Be ready also for a small-group or full-class discussion on the work that you and your classmates have done on this assignment.}
\newcommand{\handin}[2]{\emph{Before the end of class, students hand in the work they completed between Weeks\ {#1}\ and\ {#2}\ along with the cards that were used for those tasks.}}
\newcommand{\onetoone}{A single proto-guideline may be associated with only \emph{one} career, position, or job title.}

\title{HELPING STUDENTS DEAL WITH ETHICAL REASONING: THE PROTO-GUIDELINES FOR ETHICAL PRACTICE IN MATHEMATICS AS A DECK OF CARDS}

\author{Stephen M. Walk\\
Department of Mathematics and Statistics\\
St.\ Cloud State University\\
720 4th Ave. S.\\
St.\ Cloud, MN  56301\ USA\\
smwalk@stcloudstate.edu \\
(320) 308-6085
\and
Rochelle E. Tractenberg\\
Department of Neurology\\
Georgetown University\\
Suite 207 Building D\\
4000 Reservoir Rd.\ NW\\
Washington, D.C. 20057\\
rochelle.tractenberg@georgetown.edu \\
no phone number
}

\date{March 17, 2024}

\begin{titlepage}

\maketitle

\begin{abstract}
Tractenberg, Piercey, and Buell presented a list of 44 proto-Guidelines for Ethical Mathematical Practice, developed through examination of codes of ethics of adjacent disciplines and consultation with members of the mathematics community, and gave justifications for the use of these proto-Guidelines.  We propose formatting the list as a deck of 44 \emph{cards} and describe ways to use the cards in classes at any stage of the undergraduate mathematics program.  A simple game or encounter with the cards can be used exclusively as an introduction, or the cards can be used repeatedly in order to help students move to higher levels of achievement with respect to the proto-Guidelines and ethical reasoning in general; we present, in Appendix A, a sample semester-long sequence of assignments for such a purpose, with activities at various levels of Bloom's taxonomy.
\end{abstract}

\listkeywords{ethics in mathematics, ethical practice standards, proto-Guidelines, undergraduate education, capstone experiences, educational games, card games}

\section*{Acknowledgments}

This project was initiated/supported with a stipend from the NSF-funded Collaborative Research Project to Georgetown University (\#2220314), Ferris University (\#2220423), and Fitchburg State University (\#2220395). Contributions from RET were supported by the NSF-funded Collaborative Research Project to Georgetown University (\#2220314).

\vfill

\noindent Shared under CC BY-NC-ND 4.0 Creative Commons License.

\vspace{2mm}

\noindent To cite this paper: ``S.\,M.\ Walk and R.\,E.\ Tractenberg.  Helping students deal with ethical reasoning:  The proto-Guidelines for Ethical Practice in Mathematics as a deck of cards, 2024.  Preprint available at Math.arXiv (identifier and URL to be determined).''

\vspace{2mm}

\noindent The proto-Guidelines for Ethical Mathematical Practice were developed by Buell, Piercey, and Tractenberg~\cite[as cited herein]{BPT2022} and are presented in Tractenberg, Piercey, and Buell~\cite{TPB2024PG} and~\cite{TPB2024}.  They are used here by permission of that team of authors.

\vspace{2mm}

\noindent The images of the cards displayed at the end of the current paper are available as a separate .pdf file at \texttt{https://osf.io/nhq4s}.

\ 




\end{titlepage}

\section{INTRODUCTION}

\begin{quote}
\emph{``\,\ldots\,[E]thics is the effort to guide one’s conduct {\bf{\emph{with careful reasoning}}}. One cannot simply claim `X is wrong.'  Rather, one needs to claim `X is wrong because (fill in the blank).'\,''}

\hfill ---Briggle and Mitcham~\cite[p.\ 38]{BrigMitch}; \textbf{emphasis added}. 
\end{quote}

\noindent Recent years have brought an increasing interest in ensuring that mathematical practice is ethical, and in answering the following question:
\begin{quote}
\emph{1.  What does it {\small{MEAN}} to practice mathematics ethically?}
\end{quote}
For anyone whose business is \emph{teaching} those who will practice mathematics, there is a follow-up:
\begin{quote}
\emph{2.  How can the undergraduate experience help students {\small{LEARN}} to practice mathematics ethically?}
\end{quote}
Signs of interest in the first question, as well as progress toward an answer, appeared in the papers by Hersh~\cite{Hersh} and Shulman~\cite{Shulman} and the Ethics in Mathematics conferences held at Cambridge University starting in 2018.  More recently, Buell, Piercey, and Tractenberg~\cite{BPT2022} addressed the specific question of what could plausibly represent ``principles of conduct governing mathematical practice'': Ethical practice standards from the American Statistical Association~\cite{ASA} and the Association for Computing Machinery~\cite{ACM2018} were synthesized with guidance from the American Mathematical Society and the Mathematical Association of America, resulting in an initial, draft set of proto-Guidelines for Ethical Mathematical Practice that were then offered for consideration to a national sample of practitioners from the mathematics community.  The resulting set of proto-Guidelines~\cite{TPB2024}, developed after rounds of input and discussion, seem to represent the most progress so far at an answer to the first question.

Notable responses to the second question include a smattering of actual ethics in mathematics courses, whether for credit (see Webber~\cite{Webber}, for example) or not (the ``EiM'' course at the University of Cambridge is an ``informal, non-examinable course''; see~\cite{ChiodoandBursillHall}), but such courses are rare.  The proto-Guidelines just mentioned~\cite{TPB2024} can be used in the integration of discipline-specific ethics content into undergraduate courses, an example of which is the focus of this paper, so they also contribute to answering to the second question.

Additional guidance for the second question can be found in Tractenberg et al.~\cite{Tractenbergetal2017} (following the ``Mastery Rubric'' framework of Tractenberg and FitzGerald~\cite{TractFitz}) in the form of a \emph{Mastery Rubric for Ethical Reasoning}.  This Mastery Rubric involves six aspects of knowledge, skills, and abilities---the authors denote such aspects ``KSAs'' for short---in the process of ethical reasoning.  Facility with each KSA can be learned, developed, and improved over a multi-stage trajectory that is observably and concretely described using Bloom's Taxonomy of Cognitive Complexity~\cite{Bloom1956}.  Among the KSAs that constitute ethical reasoning, the first is \emph{prerequisite knowledge}\,\footnote{\ \ ``Prerequisite knowledge'' as a first aspect of reasoning is consistent with mathematics-specific learning in general; for instance, ``the knowledge base'' is the first of five parts of the mathematical-thinking framework presented in~\cite[p.\ 16]{Schoenfeld}.}; a focus on the ``prerequisite knowledge'' KSA is necessary to lay the foundation for students to engage with Briggle and Mitcham's requirement to ``guide one’s conduct with careful reasoning.''  

While some instructors may find it challenging to envision introducing prerequisite knowledge of proto-Guidelines for Ethical Mathematical Practice into their undergraduate mathematics courses, there are three contexts that are readily recognized in any such course, and which offer concrete examples of potentially unethical mathematical practice:
\begin{itemize}
\item  For convenience, conditions are assumed to hold, but they are not actually met.
\item  Approximations are used that are imprecise or simply incorrect.
\item  Familiar techniques or convenient models are used despite being inappropriate to the situation.
\end{itemize}
Questions related to such situations routinely appear in various forms---on exams, in homework, or in class sessions---for even a first- or second-semester calculus exam:  \emph{What are the hypotheses?  Verify them for this instance}, or \emph{What is the error in your approximation?  Is your answer reasonable?},\ or \emph{Is this a situation in which the model is appropriate or gives reasonable results?}  Even basic ``knowledge'' questions such as \emph{Does $3^x$ represent the same kind of function as $x^3$?  Should you be able to differentiate it via the ``power rule'' to obtain ``$x \cdot 3^{x-1}$''?} can contain traces of, say, the third consideration about appropriate models.  So the three contexts listed above can represent opportunities available in any undergraduate mathematics course for introducing aspects of ethical mathematics use and practice. Since attention to these assumptions, approximations, and applications is already paid by careful instructors and students in entry-level mathematics classes; they can be leveraged to bring ethical mathematical practice to other mathematical habits of mind and keep the inquiry
going into other classes.

Another aspect of prerequisite knowledge is an understanding of individuals who have, or who might have, an interest or stake in the outcome of mathematical work---and the certain or potential harms or costs to these individuals.  A stakeholder analysis (as introduced in~\cite{TractenbergMR}) provides a structure to consider real or potential harms or costs, \emph{and} real or potential benefits, that might accrue to a range of different people and entities---\emph{stakeholders} or, as we will call them here, \emph{interested parties}---when an ethically questionable decision is made, say, one of the three types outlined above.  Key parties to consider are the individual making the decision, their boss or client, their colleagues, the profession, and the public or public trust; further, there may be other parties, unknown at first, that are revealed through consideration of the scenario.\,\footnote{\ \ See~\cite{TractenbergLOQuantitative} for a discussion of such an analysis and learning outcomes for using it in quantitative undergraduate courses. Note also Hersh's delineation of the three broad categories of ``ethical demands of all the scientific groups'': ``what you owe the client, what you owe your profession, and what you owe the public''~\cite[p.\ 13]{Hersh}.}   
Shulman~\cite{Shulman}, for instance, cites a classic ``oil pipeline'' optimization problem, for which, though the path of the pipeline is generally chosen in order to minimize cost, there are other factors to consider---the environmental and human impact.
Consideration of harms or costs, and of the parties to whom those accrue, can help students to recognize that the mathematics they are learning is being learned in order to be applied, that the applications will involve decisions, and that decisions have consequences that are not necessarily favorable for all involved.  While the pursuit of an interested-party analysis is not currently universal to or even widespread among mathematics courses, certainly there \emph{is} a widespread desire to show students applications of the mathematics they are learning.  Thus, though the explicit inclusion of interested-party analysis may be new to any given course, it fits in with a greater theme and can be folded into inquiries that are already taking place.

These aspects of prerequisite knowledge, then, are readily available to mathematics instructors in forms already natural to or adjacent to common practices in mathematics classes.  A further aspect of prerequisite knowledge is \emph{knowledge of the ethical guidelines involved in one's field}, and because the recently-developed set of proto-Guidelines is available, instructors have a ready source to share with students.  There are obstacles:  Because there are 44 proto-Guidelines, building familiarity with them \emph{all} seems like a mammoth task, and it seems too
much to expect of students, especially on their first encounter with the list.  Worse, such an expectation might be seen by students as an exercise in memorization, and simply memorizing the set of 44 elements will not facilitate learners in learning to (recalling Briggle and Mitcham) ``guide [their] conduct with careful reasoning.''
To mitigate the difficulties posed by a consideration of the entire list, an instructor could select a specific subset of proto-Guideline elements to focus on for a single assignment or batch of assignments; a drawback to that approach is that the very act of selection risks implying that those particular proto-Guideline elements have a 
special status, as if they are the most important ones, and individual instructors might not want to risk giving that impression, or of somehow choosing the ``wrong'' ones as important.

Alternatively, students could be assigned to choose for themselves, determining their own smaller subset to focus on.  This would be an empowering, constructivist option, but it would carry the same down side noted earlier, that is, it could still be a daunting and time-consuming task, too much so for early-stage learners in their first engagement with ``ethical mathematical practice'' as a curricular objective.

As a solution to these difficulties, we describe in the current manuscript an activity that is meant to help students become familiar with the proto-Guidelines by focusing on only a subset of them, but for which the subset is determined randomly---as a deal from a deck of ``proto-Guidelines cards.''

We present, in Section~\ref{section:rulesofthegame}, details and design of the proto-Guidelines deck and basic ideas for gameplay in a single class session.  The intended student learning outcomes for a single session of basic play and a follow-up assessment are these: 
\begin{itemize}
\item  Name occupations that use mathematics explicitly or for which mathematical 
thinking is valuable.  \bloom{1}
\item  Explain how mathematics or mathematical thinking is used in those occupations.  \bloom{2}
\item  Select proto-Guidelines that might be particularly 
relevant for mathematical practice in various fields.  \bloom{1}
\item  Describe scenarios in which the selected proto-Guidelines might be particularly 
relevant for mathematical practice in various fields.  \bloom{1}
\end{itemize}
Given these objectives, one or more encounters with the deck of cards are, we believe, reasonable and positive first encounters with the proto-Guidelines.  A modest amount of preparation is required from each student---brief research into one or more careers, job titles, or employment listings---along with perhaps half a class period for the actual gameplay, followed by, if desired, a follow-up assessment such as a brief write-up, presentation, or discussion.

As noted earlier, the KSAs for ethical reasoning can be learned and developed over time.  In Section~\ref{section:variations}, we note sources that can help guide that development, and we describe variations of the game that can be applied throughout a single semester course or across various courses to help achieve that development.  The intended student learning outcomes for a larger sequence of assignments and activities include these:
\begin{itemize}
\item  Analyze scenarios (relevant to mathematics practitioners) with respect to 
possible benefits and costs to interested parties.  \bloom{4}
\item  Analyze scenarios with respect to various ethical frameworks.  \bloom{4}
\item  Recommend a course of action in scenarios where ethical decisions need to be 
made.  \bloom{6}
\item  Support action recommendations with respect to benefit-cost analysis and 
specific ethical frameworks.  \bloom{6}
\end{itemize}
In Appendix A, we present a sample sequence of assignments and class activities to illustrate some of those variations throughout a single semester course.

\section{RULES OF THE GAME \label{section:rulesofthegame}}

Our subsection headings follow the example of~\cite{Hoyle}.

\subsection{Players}

The game is played by groups of three or four.  The players can be mathematics majors in a senior capstone course, as they were in the initial proof-of-concept of the activity.  They can be mathematics majors at the beginning of their program, in an introduction/community-building course that also includes a look at jobs or careers in an attempt to help students envision what they are working toward.  They can be mathematics majors in any course in between.  They can be students of other disciplines who are taking mathematics courses as part of their program.  They can be students in a mathematics club using one or two meetings for a joint exploration of careers and ethics.  They can be, in short, \emph{any} mathematics learners who would benefit from thinking about the mathematical practice that they will be doing in their professional lives and how to carry out that mathematical practice ethically.  

Because the gameplay involves thinking about specific proto-Guidelines elements and their relevance to selected occupations and job titles, students will need to have done some research ahead of time for the activity to be meaningful.  It is essential that students have time to research the occupations and job titles, think about what they have found, and distill the research into coherent information that they share with the group.  Thus, while the \emph{players} can be any undergraduate mathematics learners, the best \emph{setting} for the game is a course or event for which careers---or, more precisely, occupations, job titles, and employment opportunities---are already an explicit part of the exploration.  (It is reasonable to have such an exploration as part of the undergraduate curriculum; the MAA's \emph{Common Vision} report, for instance, asserts that ``Our community must prepare graduates who are career-ready and focus intentionally on workplace skills early in their programs''~\cite[p.\ 20]{CommonVision}.)   For example, the\ \firstorsecond\ author has used an early version of this game in a capstone course in which roughly thirty per cent of the activities and grading are related to such exploration of employment opportunities.

Targeted assignments to explore such topics can involve sources such as the following:
\begin{itemize}
\item  A recent edition of the MAA publication \emph{101 Careers in Mathematics}~\cite{101Careers},
\item  The ``WeUseMath.org'' lists of ``Careers in Math'' and ``Careers Using Math'' at 

\texttt{http://weusemath.org/?page\underline{\ }id=143818},
\item  A search using mathematics-related search terms on \texttt{USAJobs.gov},
\item  The list of ``Dream jobs'' on the MAA Web site at  \texttt{https://mathcareers.maa.org/taxonomy/term/65} , 
\item  The SIAM ``Careers in Applied Mathematics'' Brochure at

\texttt{https://www.siam.org/Portals/0/Student\%20Programs/Thinking\%20of\%20a\%20Career/brochure.pdf}, or

\item  An online job-posting service or career-services platform with access provided by the students' own institution.
\end{itemize}
A note about terminology is order.  The various sources we suggest for students (and other sources that are not on the current list but that the instructor may suggest) present information in different ways, for different purposes, and using different language: \emph{career}, \emph{occupation}, \emph{position}, \emph{job title}, for instance.  All can be gathered under the umbrella of \emph{what students might engage in after college in order to make a living---and perhaps to make a difference}.  With this in mind, we will aim to use the related terms reasonably accurately (e.g.\ not using \emph{career} when we mean \emph{currently available job listing}), but we will avoid fine distinctions and belabored repetition of all possible terms, in an attempt to avoid being tiresome for the reader.

\subsection{Cards}

There are 44 cards, one for each of the ethical mathematics practice proto-Guidelines elements, separated into four suits according to the four categories delineated in~\cite{TPB2024}.  The separation into suits is not strictly necessary for basic gameplay,\,\footnote{\ \ The assignment of suits to the four categories of proto-Guideline elements carries meaning, though distant, related to the folklore surrounding the historical meanings of the suits, as suggested in~\cite[p.\ 41]{Hargrave} and~\cite[p.\ 62]{Tilley}.  Clubs, which represented various aspects of agriculture---loosely interpreted here as the broader world and society outside academia---have been assigned to the category of proto-Guideline elements that apply \textbf{in general} (proto-Guideline elements 1--12).  Hearts, which represented the church---loosely interpreted here as a group of people working within and toward a higher cause---have been assigned to proto-Guideline elements 13--22, which relate to the mathematics practitioner \textbf{as a member of the [mathematics] profession}.  Diamonds, which could represent the merchant classes, vassals working for a lord, or even the aristocracy, are assigned here to the proto-Guideline elements that relate specifically to projects undertaken as part of work-for-pay or making a living, that is, the ones that apply to ethical mathematics practitioners \textbf{in their scholarship} (proto-Guideline elements 23--33).  Finally, spades, representing knights or the military---notes of leadership, discipline, rank---have been assigned to those proto-Guideline elements, numbered 34--44, that apply to anyone who is a \textbf{leader/employer/supervisor/mentor/instructor} within the profession.  Other assignments could have been made to fit these historical interpretations reasonably well, but the chosen assignment has the added benefit of \emph{almost} matching the rank order of an actual card game:  It is just one transposition away from the (alphabetical) club-diamond-heart-spade ordering of bridge.} but it allows the instructor to establish house rules that encourage breadth and depth of exploration over the course of a semester.  For example, an overall goal can be stated as ``Make sure to play at least one card of each suit this semester, with at least four cards played in one particular suit,'' which is more memorable and less cumbersome than ``Make sure to play at least one card pertaining to proto-Guideline elements 1 to 12, and another pertaining to proto-Guideline elements 13 to 22, and \ldots\,.''

\begin{figure}[ht]
\centering
\includegraphics[scale=0.42]{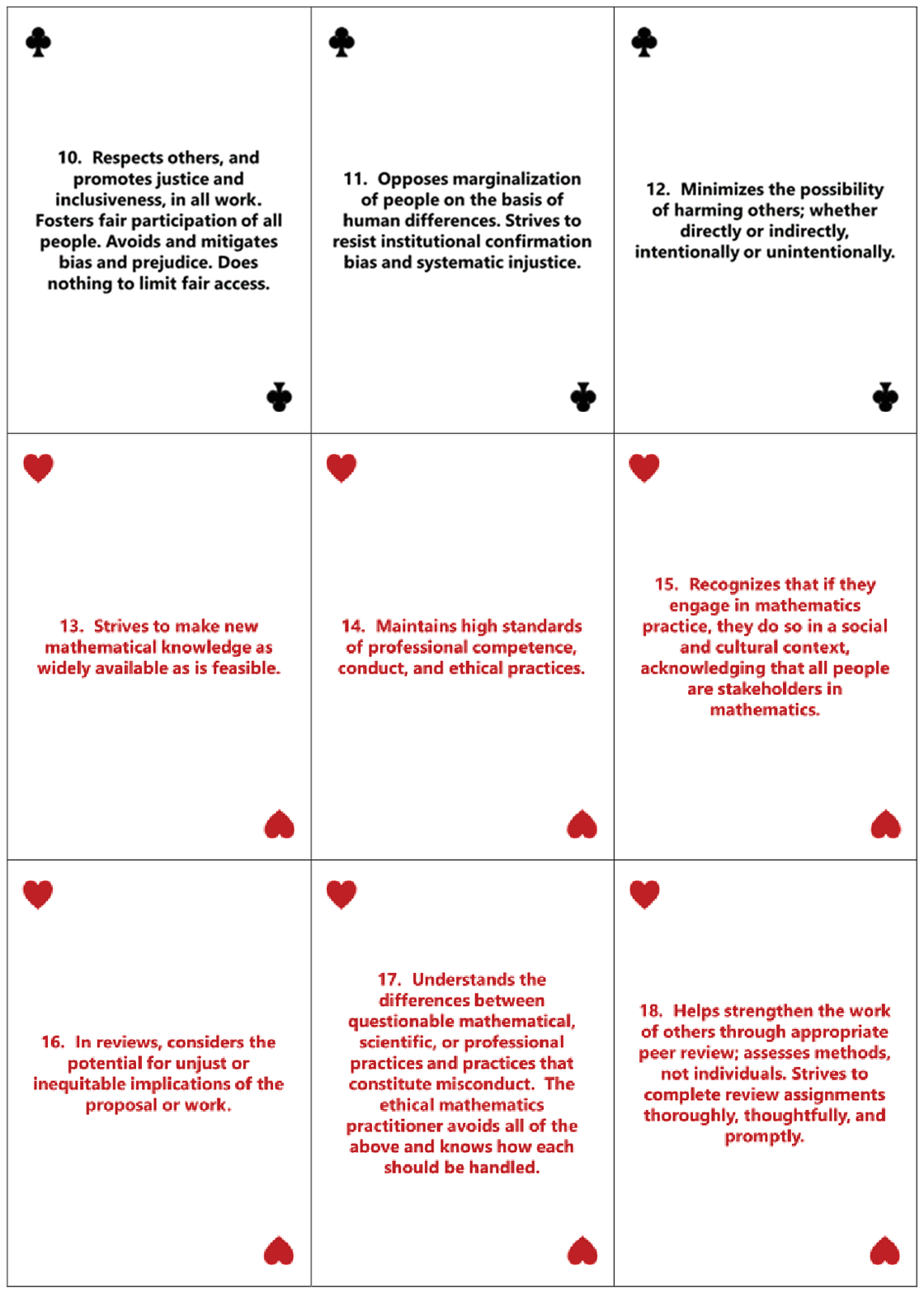}
\caption{ \label{figure:secondnine}One page of cards (nine of 44) corresponding to the second set of nine Guideline elements.  This illustration is an example; a printable .pdf version of the entire deck is in Appendix B.  The four suits of the deck correspond to the four categories of proto-Guideline elements: elements 1--12 address the behavior of the ethical mathematics practitioner ``in general," 13--22 address their behavior ``as a member of the profession,'' 23--33 address their behavior ``in their scholarship,'' and 34--44 address the \emph{additional} responsibilities of one who is ``a leader, employer, supervisor, mentor, or instructor.''  It is up to the instructor to provide this context for students.}
\end{figure}

\begin{figure}[ht]
\centering
\includegraphics[scale=0.42]{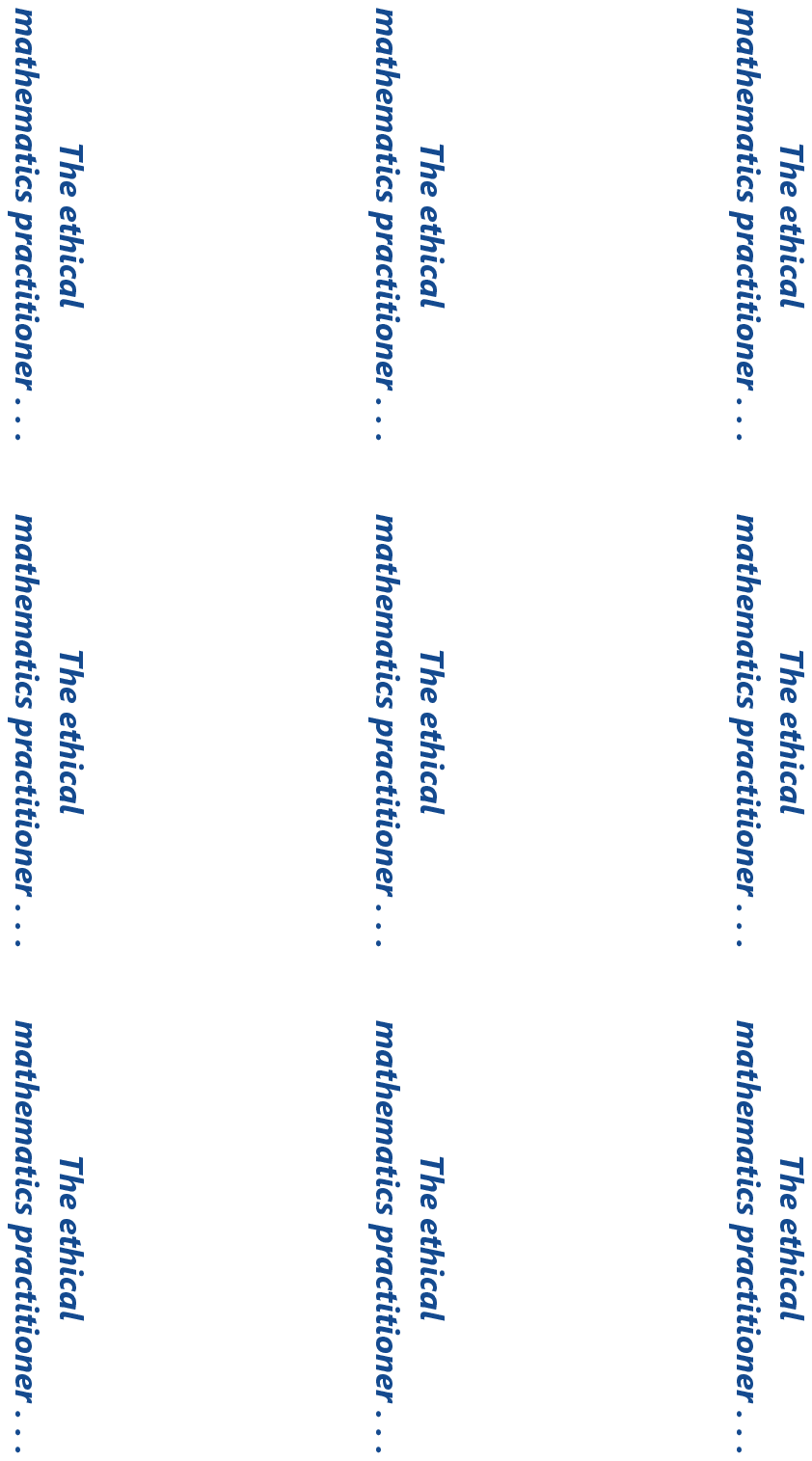}
\caption{ \label{figure:reverses}A sheet of backs for the proto-Guidelines cards.  The proto-Guideline elements for Mathematical Practice have two common stems:  Proto-Guideline elements 1--33 use the stem ``The ethical mathematics practitioner\,\ldots\,,'' while 34--44 use the stem ``An ethical mathematics practitioner who is a leader, employer, supervisor, mentor, or instructor follows all of the above items and also\,\ldots\,.''  The first, shorter, more general stem is used for all cards---as with any card game, the backs of the cards should all look the same---so it is up to the instructor to explain the context.}
\end{figure}

A printable .pdf version of the cards is included in Appendix B.  (See also \texttt{https://osf.io/x5ur9/} for the .pdf version and also the list of proto-Guidelines.)  We recommend printing the cards on white cardstock of density at least 176 grams per square centimeter.  Each sheet needs to be printed on twice---once for the card faces and once for the card backs---and then there is cutting to do, as each sheet displays eight or nine of the 44 cards.  Color printing is nice but not necessary.  Figure~\ref{figure:secondnine} shows one page of nine faces; Figure~\ref{figure:reverses} shows the page of card backs.

\ 

\subsection{The Deal}

Each player is dealt seven cards at the beginning of play.  For later hands (later in the semester, that is), the instructor will probably need to include one or two additional decks per group to ensure enough cards to go around; if that step is taken, then rules about repeated cards should be implemented (e.g.\ ``If your hand contains two of the same card, you must discard one and receive a new card'' and ``You may not play two copies of the same card during the semester'').

\subsection{The Play}

As noted, each student comes to class having done some research on a certain number of careers.  The basic idea of play is that the student matches each of their researched descriptions of employment with a proto-Guideline element card from their hand.  It is not absolutely necessary to take part in a class activity in order to do this matching; each student can look at their own hand of cards and select the ones that seem most relevant, and this can be done outside of class if preferred.

A way to make it more fun and interactive---and to encourage every student to think about more proto-Guideline elements than just the seven in their hand and more employment possibilities than just the ones they did their own research on---is to use an approach similar to that of numerous well-known party games:
\begin{itemize}
\item  Player A, the student chosen to go first, summarizes or recaps, for all 
players in their group, one of their researched careers or occupations or job titles, describing, for 
instance, its duties, its required or recommended skills, and some other aspect that
Player A found interesting.
\item  Each of the other players chooses from their hand one ethical proto-Guideline element
that someone in that position might need to be especially aware of. The chosen cards are 
placed face down in front of Player A.
\item  After all players have played their cards, Player A reads the played cards 
aloud and chooses one, laying it face up and giving a brief rationale for choosing 
it over the others. Player A puts that card face up into their own ``chosen 
cards'' pile; the other played cards are shuffled back into the larger pack.
\item  Players pick cards to replenish their hands, and the turn passes to the next 
player.
\item  Play continues for a number of turns determined by the instructor.
\end{itemize}
The game has no defined ``winner'' but allows all students to benefit by engaging 
with each other and with the proto-Guidelines. Students leave class with a full 
hand of seven cards \emph{plus} the cards from their ``chosen cards'' pile. 
They are assigned to complete the following steps outside of class (details of a sample assignment are given in Appendix A):
\begin{itemize}
\item  For each of the researched employment descriptions, decide which of the currently-held cards is really the best fit---the proto-Guideline element that seems the most relevant.  It does not have to be one of the cards that were 
played and chosen during the game; it can be one of the cards that the player already held in their hand.
\item  Complete some sort of assessment---for instance, a write-up of a few 
sentences---for each of the career-card pairings just created.
\end{itemize}
It is a good idea for the instructor to let students in 
on one piece of strategy: Since they will ultimately choose fitting 
proto-Guideline elements from among \emph{all} of the cards they take from the classroom, it is wise \emph{not} to play 
a card from their hand in class if they perceive immediately that it is a really good fit for one of 
the positions or careers that they researched.

\section{VARIATIONS \label{section:variations}}

\begin{quote}
\emph{…[E]thics is not a vaccine that can be administered in one dose and have long lasting effects no matter how often, or in what conditions, the subject is exposed to the disease agent.}

\hfill ---National Academy of Engineering~\cite[p.\ 34]{NAE2009}
\end{quote}

\noindent As described above, the card game can be played in just a single class session, or
just a single mathematics club meeting, as an introduction to the proto-Guidelines. It can
be played exactly the same way, repeatedly, over the course of a semester as a means to keep the proto-Guidelines at the front of students’ minds.

Even with a set of accepted ethical practice guidelines, integrating ethical issues throughout a curriculum requires a focus on careful reasoning with those guidelines; we recall Briggle and Mitcham's quote that opened the Introduction, ``Ethics is the effort to guide one’s conduct with careful reasoning.''  We turn now to ideas for repeated work with the proto-Guidelines deck that will help students consider the proto-Guidelines more deeply and develop their skill at ethical reasoning.  Sources to guide these efforts include the Value Rubric for Ethical Reasoning published by the American Association of Colleges \& Universities (see for example~\cite[pp. 46--47]{ValueRubrics}), which lists four levels from ``Benchmark'' to ``Capstone,'' and the aforementioned Mastery Rubric for Ethical Reasoning proposed by Tractenberg and FitzGerald~\cite{TractFitz} and elaborated by Tractenberg, FitzGerald, and Collman~\cite{Tractenbergetal2017}, with levels from ``Novice'' to ``Master.''  While the top level of the AAC\&U VALUE Rubric, ``Capstone,'' is meant to be achieved within a student's undergraduate program, the top levels of the Mastery Rubric, ``Journeyman'' (independent practitioner) and ``Master'' (qualified to take on an apprentice), are meant to apply to those already working in the field, with an undergraduate education bringing the practitioner only to the level of ``Apprentice'' in anticipation of entering the workforce.  The levels of ``ethical engagement'' described by Chiodo and Bursill-Hall~\cite{ChiodoBursillHallFourLevels} are worth noting---four levels \emph{plus} a ``Level 0'':  Level 1, ``Realising there are ethical issues inherent in mathematics,'' is the highest to expect for an undergraduate experience, though Level 2, ``Doing something: speaking out to other mathematicians,'' is not entirely out of reach for undergraduates on the bold side.

Also relevant is the idea of \emph{stewardship of a discipline}:  A \emph{steward}, as defined by Golde and Walker in their edited volume exploring the construct~\cite{GoldeWalker}, is ``an individual to whom we can entrust the vigor, quality, and integrity of the field.''  They note further that
\begin{quote}\emph{\ldots\,upon entry into practice, all professionals assume at least a tacit responsibility for the quality and integrity of their own work and that of colleagues. They also take on a responsibility to the larger public for the standards of practice associated with the profession.}~\cite[p.\ 10]{GoldeWalker}\end{quote}
As noted by Tractenberg in~\cite{Tractenberg2022}, instructors who are hesitant to structure \emph{ethical content} for mathematics curricula could instead consider encouraging students of mathematics to be \emph{stewardly} in their use of the tools and techniques of the discipline.  The description of a steward in~\cite[p.\ 5]{GoldeWalker} was originally articulated for doctoral students, but Rios et al.~\cite{RiosetalStewardship} expanded the definition to include students at all stages of professional development, and an emphasis on stewardly use of mathematics has clear relevance for mathematics instruction that seeks to prepare practitioners of mathematics---or \emph{any} users of mathematics.\,\footnote{\ \ In his essay ``What Is Mathematics For?''~\cite{DudleyWhatIsMathematicsFor},\ Underwood Dudley refers to the concept of a \emph{steward of mathematics}---briefly, but to great effect.}

Building on the Mastery Rubric for Stewardship presented in~\cite{RiosetalStewardship}, Tractenberg~\cite{Tractenberg2022} articulated learning objectives specific to four of the eight KSAs of that Mastery Rubric, noting that ``four of the stewardship KSAs can be brought to bear throughout \textbf{most} mathematics courses'' (emphasis added).  Reworded slightly from~\cite{RiosetalStewardship}, these four are the following: 
\begin{itemize}
\item  Requisite knowledge and situational awareness, 
\item  Critical evaluation of extant knowledge,
\item  Responsible application of disciplinary knowledge, and 
\item  Responsible communication.
\end{itemize}
Thus, an instructor or team of instructors who emphasize these aspects of problem-solving in their courses are helping their students to become stewards of the discipline.

In an entry in the Golde and Walker volume~\cite{GoldeWalker}, Bass~\cite{Bass2006} considered the construct of stewardship within the training of scholars and professionals in the field of mathematics specifically, asking, ``In what ways does a mathematician function as a representative of the discipline in public arenas?''~\cite[p.\ 111]{Bass2006}\ \ The instructor interested in leveraging stewardship will easily perceive the relevance of this question for the user of mathematics at work, whether or not their job title is ``mathematician,'' and hence for undergraduates along their instructional trajectory, whether or not they continue on to graduate study.

Following curriculum development guidelines such as in~\cite{Nicholls2002} and~\cite{TractenbergMR}, three learning outcomes might be that, after a course with explicit stewardship and proto-Guidelines content, the students will be able to do the following:
\begin{enumerate}
\item  \emph{Describe how a given mathematics course contributes to their sense of professional identity, and compare the contribution of such a course to this sense (with and without the construct of stewardship).}
\item  \emph{Relate stewardship KSAs to the conceptual content of the course---for example, proofs, definitions, theorems, applications.} 
\item  \emph{Justify how many mathematics courses must be taken before a student owes a duty of stewardship to the discipline of mathematics.}
\end{enumerate}
The Mastery Rubric for Stewardship itself has four stages, each concretely and distinctly defined according to the level of complexity and sophistication with which individuals perform each KSA. The earliest level in this Mastery Rubric, ``Novice,'' is described as ``an early stage learner who does not recognize that, or act as if, failures to act in a stewardly manner have ramifications beyond themselves. The novice stage represents the individual embedded in the acquisition of discipline-specific content'' (see~\cite{RiosetalStewardship}; with an eye toward ethical reasoning in particular, this stage corresponds to the Level 0 of Chiodo and Bursill-Hall~\cite{ChiodoBursillHallFourLevels}, whereas the lowest level of the AAC\&U VALUE Rubric, ``Benchmark,'' assumes that the student can already recognize ethical issues).  The next stage of performance, ``Apprentice,'' is described as ``an individual who is actively engaged in study of the profession or discipline\,\ldots\,developing the capacity to practice independently, but has not yet demonstrated ability qualification to do so.''

It was noted earlier that one author of the current paper used the basic version of gameplay with the cards in a session of a capstone course for senior mathematics majors.  The follow-up assessment was a brief write-up in which students explained why they believed that certain proto-Guidelines were relevant to their researched occupations.  That author was pleased enough with the results that he intends to use the activity again, but in an expanded form, as a sequence of skill-building activities with multiple variations, undertaken over the course of a semester in an attempt to help students move from the ``Novice'' level to the ``Apprentice'' level (of stewardship~\cite{RiosetalStewardship} or of ethical reasoning~\cite{TractFitz}), or beyond the ``Benchmark'' level and to the ``Capstone'' level of the VALUE Rubric for Ethical Reasoning (see~\cite[pp.\ 46--47]{ValueRubrics}).  Some considerations for repeated, scaffolded gameplay are as follows:
\begin{itemize}
\item  To help prevent student overload and decision fatigue, the sources recommended 
for career and employment research can be limited at the beginning of the semester and expanded as 
the semester goes on---for example, students might be restricted to look in just one or two 
specific sources the first time around, the list could be expanded to four sources 
the next time around, and perhaps students could seek out their own sources the third 
time as desired.
\item  An overall goal can be given for the semester: say, ten cards played in all, 
with a requirement such as ``you need a flush---cards all of the same suit---of a 
certain size \emph{and also} make sure you play cards of all four suits.'' A gameboard or 
scoreboard can be created by the instructor---or by the students!---to help students keep track.
\item  As students continue to explore more career options and more sources, they can build on 
previous explorations and dig more deeply into work done earlier 
in the semester. For instance, students can be assigned to do the following for one
or more of their explored careers, with each bulleted item building on work done in 
previous bulleted items as directed by the instructor: \em
\begin{itemize}
\item  Brainstorm and fill out an interested-party analysis for a course of action 
related to the scenario.
\item  Pick two different ethical frameworks: virtue ethics, utilitarianism, natural 
law, social contract, for instance. (Some possible frameworks will be outlined in 
class, with examples given, and a written overview of those will be provided. 
Students who are aware of other frameworks, and who wish to use those, are welcome 
to do so.) As before, brainstorm. For each chosen ethical framework, use the 
interested-party analysis to imagine how the proposed action would be viewed by a 
proponent of that framework.
\item  Answer the following questions: Are there alternative actions [beyond just 
``do or don’t do the action that was proposed'']? Which of the interested parties 
might consider each such action most relevant or compelling? What would adherents 
of two or three different ethical frameworks say? 
\item  Propose a course of action that seems reasonable for an ethical mathematics 
practitioner in the situation. Describe how the mathematics practitioner could 
defend the decision, specifically to those with other ethical frameworks or to those
who incur costs from the proposed course of action. Then tell how and when the 
mathematics practitioner might reasonably reflect on the course of action---considering the 
ramifications of the action and implications for similar decisions in the future---and 
criteria that might be used.
\end{itemize}
\end{itemize}
Moreover, in anticipation of a final exam or final project that includes engagement with the Guidelines, ethical reasoning, a stewardly approach to mathematical practice, or any combination of these, an instructor can ask for short essays (``minute papers'') throughout the term wherein students reflect on the games that were played and any work that they produced that pertained to the games.  Prompts could be, for example, ``Describe two things you learned (in this game) that were unexpected, and tell why you were not expecting them" or ``List two occupations that were not considered (today) where mathematics might be important, and briefly outline why the practitioner of mathematics in those occupations should or should not follow the same Guidelines for Ethical Mathematical Practice as someone with an occupation that featured mathematics more prominently.''  With a term's worth of such essays completed, students can be asked to reflect on the essays, describing how their thinking about the idea of ``ethical mathematical practice'' has changed, as reflected in their essays (i.e.,\ with evidence they glean from their own work); or to discuss the role the card games played in the mathematics work they completed for the course; or to reflect on the \emph{relevance} versus \emph{prevalence} of these ideas throughout their other mathematics courses to date.

A sequence of sample assignments, representing such a deeper dive over the course of a semester, is provided in Appendix A.  Alternatively, activities using these variations can be undertaken in a selection of courses throughout the mathematics major program rather than in a single semester.

As noted earlier, students do not need to \emph{play a game in class} in order to use the cards effectively; each individual student can select cards on their own, outside of class if preferred, and complete the assignments just as if they had received extra cards from classmates by playing the game in class.  They can bring their written assessments
to class, receive one or more cards to replace the card(s) just used for those 
assessments, and continue to make selections on their own rather than actually 
playing the game. (In keeping with the rationale of the game, it is best not to 
overload students by requiring them---or even allowing them---to pick from the entire list but, rather, to require them to pick 
from the smaller number of cards randomly dealt to them or contributed by fellow 
players.)

\section{CONCLUSION}

Given the ubiquity of mathematics throughout scientific and technological jobs and applications globally, it is vital for those workers in those fields to practice mathematics ethically \emph{whenever they practice mathematics}, even if that is not their primary focus, just as it is vital for them to be ethical in their use of technology when they use technology, ethical in their use of data when they work with data, and ethical in their communication when they communicate, regardless of \emph{how much} of their working time is spent using technology, analyzing data, or communicating.  To that end, it is worthwhile for students to consider the proto-Guidelines of Ethical Mathematical Practice with respect to a wide range of fields in which mathematics \emph{can be} used, even if developing or applying new mathematics is not the \emph{primary focus} of such fields.  For those whose work \emph{is} primarily about developing or applying new mathematics, the following is relevant:
\begin{quote}
\emph{The power of new mathematics in ethically-laden industries means the professional and temporal gap between its creation and its application has reduced so much that the ethical consequences of mathematical work cannot be obscured or blamed on someone else.  For {\bf the first time ever}, mathematicians are {\bf uniquely} responsible for the immediate social consequences of their work.}

\ \hfill ---Chiodo and Bursill-Hall~\cite[p.\ 22, emphasis in original \emph{as italics}]{ChiodoandBursillHall}
\end{quote}
In the same essay, Chiodo and Bursill-Hall note the need for students to ``train their ethical reasoning just like they train their mathematical reasoning via exercises''~\cite[p.\ 25]{ChiodoandBursillHall}.  We cannot expect that students  who have never practiced or even seen ethical reasoning to suddenly become masters at ethical reasoning when they enter the workforce.  

This is not to say that it is necessarily an easy road to introduce students to ethical reasoning within a mathematics program or a mathematics course---or to win colleagues over to the idea.  It might be that some colleagues are at Level 0 of awareness of ethical concerns: Chiodo and Bursill-Hall remark that ``One reason mathematicians shy away from ethical discussions is that mathematics seeks timeless, absolute truths.\,\ldots\,[E]thics doesn't have the same binary clarity or timelessness''~\cite[p.\ 25]{ChiodoandBursillHall}.  Another possibility is that because the disciplinary content for the mathematics undergraduate curriculum is already so substantial---some might say \emph{overwhelming}---that adding ethical reasoning to the mix, even if such reasoning is focused on discipline-specific guidelines, is challenging; asking instructors to include ethical content might suggest that some of the disciplinary content should be removed to make way for it. 

Our intent in describing the proto-Guidelines deck, and the ideas for activities and assessments that use it, is to provide a means for students to build awareness of the proto-Guidelines with minimal use of class time (and by leveraging instruction that is already, typically, in use), or to develop a higher level of ethical reasoning ability with repeated engagement, at whatever level is palatable for the individual instructor and realistic for the specific program.    Experience with the proto-Guidelines is valuable; as instructors, we just have to have a sense of the best times, within the scope of a program, to hold them up as worthy principles.  Moreover, we believe that activities that provide practice in ethical reasoning, and engagement with the proto-Guidelines for Ethical Mathematical Practice, can be reasonably pursued within one or more existing mathematics courses; we need only to be alert for opportunities to fold them into the curriculum.  Then we can see our students as future mathematics practitioners, we can raise their awareness of ethical responsibilites, and---ideally---we can begin to call them stewards of the discipline.

\vfill

\ 

\pagebreak

\section*{APPENDIX A}

We present a sample sequence of actvities and assignments that can be carried out over the course of a semester.  The assignments incorporate many of the variations suggested in Section~\ref{section:variations}.  Deepening exploration can be achieved by using the same sorts of tasks, but asking students to perform at increasingly higher levels of Bloom's taxonomy~\cite{Bloom1956} or using the Mastery Rubric for Stewardship in~\cite{RiosetalStewardship} (see also~\cite{Tractenberg2022}).  Another way to deepen student engagement with ethical considerations of using mathematics practices in the workplace is to use the Ethical Reasoning VALUE Rubric in~\cite{ValueRubrics}.

The tasks representing the increasingly higher levels of performance are labeled ``Rank 1'' to ``Rank 4'' in the sample assignments; this ranking is internal to the sequence of assignments and does not correspond directly to Bloom's taxonomy or to any of the rubrics mentioned earlier.  The overall goal for students is to ``play'' ten cards over the course of a semester, at least one of each suit and at least four in one particular suit.  This requires them to attend to all four general sections of the Guidelines, and it facilitates discussion about why these sections exist.

The sample materials are presented in the context of a class that meets once per week.  The first assignment of the sequence, given during Week 1 of the course, is as follows.

\vspace{4mm}

\noindent \hrulefill

\vspace{2mm}

\assignment{1}{1}{2}

\vspace{2mm}

\fordeadline{one week}
\begin{enumerate}
\item  Read at least two of the profiles in \emph{101 Careers in Mathematics} and fill out the attached information sheet for each.  At least one of the two you select should be about a professional with only a bachelor’s degree in mathematics no other fields.
\item  Read at least two of the descriptions of ``Careers in Math'' and ``Careers Using Math'' on the WeUseMath.org Web site at \texttt{http://weusemath.org/?page\underline{\ }id=143818} .  Fill in the information sheet (attached) for the two that you chose.
\item  Prepare a one-minute talk (notes and information sheet allowed) about one of the people you read about in \emph{101 Careers}.
\item  Prepare another one-minute talk (notes and information sheet allowed) about one of the careers you read about on WeUseMath.org.
\end{enumerate}
Notify your instructor via e-mail about the four you selected at least 24 hours ahead of next week's class.

\vspace{2mm}

\vfill

\ 

\pagebreak

\noindent \hrulefill

\vspace{3mm}

\noindent \textbf{Information sheet for a professional profile from \emph{101 Careers in Mathematics}}

\vspace{6mm}

\noindent Name:

\vspace{6mm}

\noindent Education (degrees and fields, not institutions):

\vspace{6mm}

\noindent Position title, and company or institution:

\vspace{6mm}

\noindent Duties:

\vspace{6mm}

\noindent How mathematics is used, including (perhaps) the particular mathematics topics or skills involved in any given task for someone in this position:

\vspace{6mm}

\noindent Aspects of the position (or of the field in general) that you think you would enjoy or would find interesting:

\vspace{6mm}

\noindent Something else interesting about this position or this line of work:

\vspace{6mm}

\noindent \hrulefill

\vspace{3mm}

\noindent \textbf{Information sheet for a career description from WeUseMath.org}

\vspace{6mm}

\noindent Career title:

\vspace{6mm}

\noindent Education suggested:

\vspace{6mm}

\noindent Duties (perhaps including potential employers):

\vspace{6mm}

\noindent How mathematics is used:

\vspace{6mm}

\noindent Mathematics topics or skills that might be used in any given task for someone in this career:

\vspace{6mm}

\noindent Aspects of the position (or of that field in general) that you think you would enjoy or would find interesting:

\vspace{6mm}

\noindent Something else interesting about this career:

\vspace{6mm}

\vfill

\ 

\pagebreak

\noindent \hrulefill

\vspace{4mm}

\noindent \textbf{{\small{SAMPLE}} class activity for Week 2.}  \emph{In class during Week 2, students deliver their brief presentations.  The instructor may bring up the subject of ethical mathematical practice in general and give one or two examples of ethical issues and questions that might confront people in certain positions or occupations.  (Ideally, the instructor will choose example positions or occupations not among those presented by students; that was the reason for the 24-hour notice.)}

\emph{Students separate into groups of three or four.  The instructor provides background on the development of the proto-Guidelines and credits the Tractenberg-Piercey-Buell article~\cite{TPB2024}.  The instructor then hands a pack of proto-Guidelines cards to one student in each group and cites the current article.  The students with card packs deal out the first hand of the game, and the game is played within each group until each student has had at least two turns to present careers and choose cards.}

\noindent \hrulefill

\vspace{2mm}

\assignment{2}{2}{3}

\vspace{2mm}

\fordeadline{one week}
\begin{itemize}
\item  For \emph{\small{TWO}} of the career profiles for which you played the game in class, do the following:
\begin{itemize}
\item  Select, from among the cards ``won'' during our in-class game \emph{and} the cards still in your hand, the ethical proto-guideline that you think is most relevant for someone in the selected career or position.\ \ \onetoone
\item  For each of the careers or positions, brainstorm ways in which a person with that career or position might need to apply the proto-guideline.  Complete the ``Rank-1'' tasks on the attached worksheet and bring it to class next week with your career-description worksheet and the relevant proto-Guidelines card.
\end{itemize}
\item  Come to class ready to give a one-minute presentation (with notes and the worksheet) on each of your career-card pairings and the associated interested-party analysis.\ \ \discuss
\item  To prepare to play the game in class again, do the following:
\begin{itemize}
\item  Pick two new careers or job titles, at least one of which comes from a search 
using mathematics-related search terms on \texttt{USAJobs.gov}.
\item  Fill out an information sheet for each of them.\,\footnote{\ \ Now that examples of information sheets have been provided, we will refrain from illustrating them for every possible source.}
\item  Come to class ready to give a one-minute presentation (with notes and the worksheet) on each of the careers you selected.\ \ \discuss
\end{itemize}
\end{itemize}

\vspace{4mm}

\noindent \hrulefill


\begin{center}
\textbf{Rank-1 tasks (to be completed as homework)}
\end{center}


\noindent \emph{Career or position:}

\vspace{4mm}

\noindent State, and give the number of, a proto-Guideline element, from the cards currently in your hand, that you have identified as important to someone in that career or position.  Write two to three sentences about {\small{WHY}} you think that that proto-guideline is important for a person in the selected career or position.

\vspace{4mm}

\noindent Brainstorm and then describe, in one to three sentences, a scenario in which someone in that career or position might need to apply that proto-Guideline element to their actions.

\vspace{6mm}

\noindent Describe a decision that they might have to make, in such a scenario, that could be aided by the guidance given in the proto-Guideline element.

\vspace{6mm}

\noindent Imagine that someone in the career or position is considering making that decision a certain way, that is, the person is considering a particular course of action to pursue instead of another one.  Fill out an interested-party analysis for the situation where the individual with the identified position fails to follow the four Guideline elements you identified previously.
\begin{itemize}
    \item  Who could \emph{benefit} from that decision?  List some parties who could benefit from such a decision, telling how each could be better off.  (You may find it helpful to separate some of these thoughts into ``short-term'' and ``long-term'' benefits.)

\begin{center}
    \begin{tabular}{|c|c|} \hline
    \textbf{Who benefits?} \hspace{15mm}  & \textbf{How?} \hspace{60mm} \\ \hline \hline
         &  \\ 
         &  \\
         &  \\
         &  \\
         &  \\
         &  \\
         &  \\
         &  \\ \hline
    \end{tabular}
\end{center}

\vspace{6mm}

\item  Who could \emph{incur costs} from that decision?  List some parties who could incur costs from such a decision, telling how each could be worse off.  (You may find it helpful to separate some of these thoughts into ``short-term'' and ``long-term'' benefits.)

\begin{center}
    \begin{tabular}{|c|c|} \hline
    \textbf{Who incurs a cost?} \hspace{6mm}  & \textbf{What is the cost?} \hspace{36mm} \\ \hline \hline
         &  \\ 
         &  \\
         &  \\
         &  \\
         &  \\
         &  \\
         &  \\
         &  \\ \hline
    \end{tabular}
\end{center}
     \end{itemize}

\vspace{4mm}

\pagebreak

\noindent \hrulefill

\vspace{4mm}

\noindent \textbf{{\small{SAMPLE}} class activity for Week 3.}  \emph{In class during Week 3, students deliver their brief presentations.  The instructor may facilitate a discussion with the class.  Points to touch on during the discussion---perhaps phrased as questions to the group---could be the following:
\begin{itemize}
\item  Following the Guidelines can add extra work, or might add extra time to get the same task done, while failing to follow one or more Guidelines may save time. 
\item  An individual who follows the Guidelines strengthens the profession, while one who does not will weaken the profession.
\item  If a company or organization claims that it, and all its employees, always follow the Guidelines, it could boost that organization's reputation with the public, just as either falsely making this claim or failing to acknowledge that proto Guidelines for Ethical Mathematical Practice exist could weaken or undermine the public trust in the organization.
\end{itemize} 
The instructor may bring a new deck of cards for each group of students, shuffling the previous game's unused cards into the larger deck.  Each student receives new cards to replenish their hand from the cards used for the homework.  The game is played within each group, for the new careers selected between Week 2 and Week 3, until each student has had at least two turns to present careers and choose cards.}

\emph{The instructor provides brief introductory readings on various ethical decision-making frameworks and outlines the basic ideas of each in class.}

\handin{2}{3}

\vspace{4mm}

\noindent \hrulefill

\vspace{2mm}

\assignment{3}{3}{5}

\vspace{2mm}

\begin{figure}[h]
\centering
\includegraphics[scale=0.68]{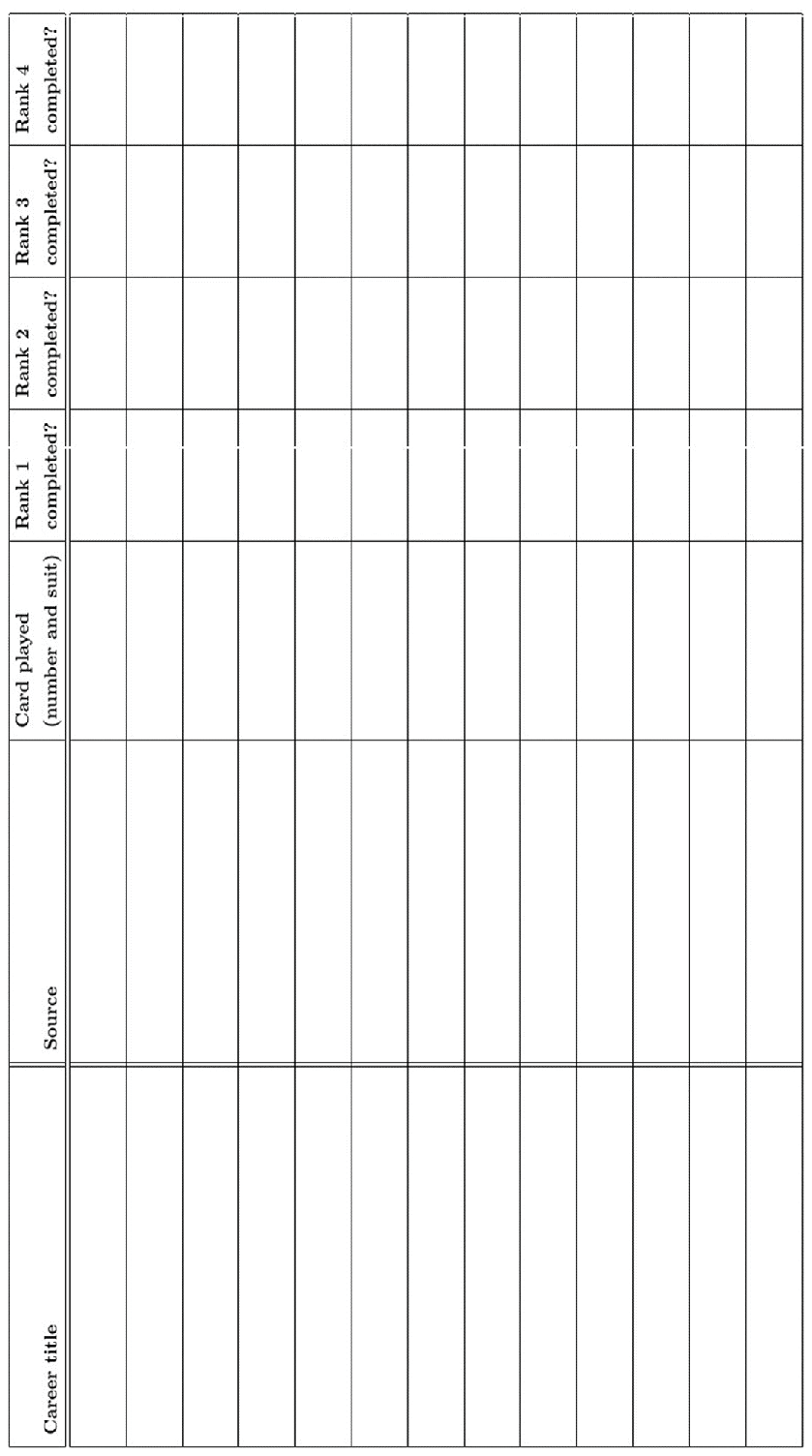}
\caption{ \label{figure:bigpictureworksheet}A ``big-picture worksheet'' for keeping track of cards played.}
\end{figure}

\fordeadline{two weeks}
\begin{itemize}
\item  For \emph{\small{TWO}} of the career profiles for which you have played the game in class but have not yet performed the Rank-1 tasks, do the following:
\begin{itemize}
\item  Select, from among the cards ``won'' during our in-class game \emph{and} the cards still in your hand, the ethical proto-guideline that you think is most relevant for someone in the selected career or position.\ \ \onetoone
\item  For those cards and the matching careers---those ``career-card pairings''---perform the Rank-1 tasks.
\end{itemize}
By the end of the semester, you are to have created ten career-card pairings in all, using at least one card from each suit and using at least four cards from one particular suit of your choice.  As you make your selections, keep that goal in mind.
\item  For two of the career-card pairings for which you have performed the Rank-1 tasks, also complete the ``Rank-2'' tasks on the attached worksheet.
\item  Bring all worksheets and relevant proto-Guidelines cards to class in two weeks.  Come to class ready to give a one-minute presentation (with notes and the worksheets) on each of your career-card pairings and the associated Rank-2 tasks. \ \ \discuss
\item  To prepare to play the game in class again in two weeks, do the following:
\begin{itemize}
\item  Pick two new careers, job titles, or occupations.  Try to include at least one that comes from a current opening posted on the online job-posting service or career-services platform connected with our institution.  (If you object to the platform's Terms of Use or Privacy Policy or object to being asked to sign up for it for any other reason, you are welcome to omit this recommendation, but do still pick two new careers.)
\item  Fill out an information sheet for each.
\item  Come to class ready to give a one-minute presentation (with notes and worksheets) on each of the new careers you selected.
\end{itemize}
\end{itemize}
\emph{Summary:}  When you come to class two weeks from now, you will have
\begin{itemize}
\item  two career-card pairings analyzed all the way to Rank-2 tasks,
\item  two other career-card pairings analyzed with respect to Rank-1 tasks, and
\item  four careers investigated but not yet matched with proto-Guidelines cards.
\end{itemize}
A big-picture worksheet is attached to help you keep track of your work.  You will be asked to hand this in at the end of the semester.  The worksheet has spaces for more than ten options in case you decide to drop your exploration of some careers or occupations in favor of others.  (The worksheet appears in this paper as Figure~\ref{figure:bigpictureworksheet}.)

\vspace{4mm}

\noindent \hrulefill

\vspace{4mm}

\begin{center}
\textbf{Rank-2 tasks (to be completed as homework)}
\end{center}

\vspace{4mm}

\noindent \emph{Career or position:}

\vspace{6mm}

\noindent Do the assigned reading on ethical frameworks (virtue ethics, utilitarianism, natural 
law, social contract, for instance) and review what was noted about them in class.  Pick two different ethical frameworks from among these.  (If you are aware of other frameworks and wish to use those, then you are welcome 
to do so.)

\vspace{2mm}

\noindent As before, brainstorm. For each chosen ethical framework, refer to the 
interested-party analysis from your Rank-1 work and imagine how the proposed action would be viewed by a 
proponent of that framework.

\vspace{6mm}

\noindent Summarize your reflections in three to five sentences per ethical framework per career or position.

\vspace{6mm}

\vfill

\ 

\pagebreak

\noindent \hrulefill

\vfill

\ 

\pagebreak

\noindent \hrulefill

\vspace{4mm}

\noindent \textbf{{\small{SAMPLE}} class activity for Week 5.}  \emph{In class during Week 5, students deliver their brief presentations.  The instructor supplies new decks of cards as needed.  The game is played within each group, for the new careers selected between Week 3 and Week 5, until each student has had at least two turns to present careers and choose cards.}

\handin{3}{5}

\vspace{4mm}

\noindent \hrulefill

\vspace{2mm}

\assignment{4}{5}{8}

\vspace{2mm}

\fordeadline{three weeks}
\begin{itemize}
\item  For \emph{\small{TWO}} of the careers or positions or job titles for which you have played the game in class but have not yet performed the Rank-1 tasks, do the following:
\begin{itemize}
\item  Select, from among the cards ``won'' during our in-class game \emph{and} the cards still in your hand, the ethical proto-Guideline element that you think is most relevant for someone in each selected career or position.\ \ \onetoone
\item  For those career-card pairings, perform the Rank-1 tasks.
\end{itemize}
\item  For two of the career-card pairings for which you have performed the Rank-1 tasks but not the Rank-2 tasks, complete the Rank-2 tasks.
\item  For two of the career-card pairings for which you have performed the Rank-2 tasks, also complete the ``Rank-3'' tasks on the attached worksheet.
\item  Bring all worksheets and relevant proto-Guidelines cards to class in two weeks.  Come to class ready to give a one-minute presentation (with notes and the worksheet) on each of your career-card pairings and the associated Rank-3 tasks. \ \ \discuss
\item  To prepare to play the game in class again in three weeks, do the following:
\begin{itemize}
\item  Pick two new careers or occupations.  These can come from any of the four sources recommended so far, or from a new source (including a personal connection or job-shadowing contact) as long as you properly credit the source.
\item  Fill out an information sheet for each.
\item  Come to class ready to give a one-minute presentation (with notes) on each of the new careers you selected.
\end{itemize}
\end{itemize}
\emph{Summary:}  When you come to class three weeks from now, you will have
\begin{itemize}
\item  two career-card pairings analyzed all the way to Rank-3 tasks,
\item  two career-card pairings analyzed all the way to Rank-2 tasks,
\item  two other career-card pairings analyzed with respect to Rank-1 tasks, and
\item  four careers investigated but not yet matched with proto-Guidelines cards.
\end{itemize}

\vspace{4mm}

\vfill

\ 

\pagebreak

\noindent \hrulefill

\vspace{4mm}

\begin{center}
\textbf{Rank-3 tasks (to be completed as homework)}
\end{center}

\vspace{4mm}

\noindent \emph{Career or position:}

\vspace{6mm}

\noindent Refer to the brainstorming, proposed actions, and possible responses from various ethical viewpoints that you completed for your Rank-1 and Rank-2 tasks.  Respond to the following, giving each response in one to three sentences.

\vspace{6mm}

\noindent What are some alternative actions to consider (besides just ``do or don’t do the action 
that was proposed'')?

\vspace{6mm}

\noindent Which of the interested parties might consider each such alternative action(s) most relevant or
compelling?

\vspace{6mm}

\noindent What would adherent of two different ethical frameworks say about each alternative action? 

\vspace{6mm}

\noindent \hrulefill

\vspace{4mm}

\noindent \textbf{{\small{SAMPLE}} class activity for Week 8.}  \emph{In class during Week 8, students deliver their brief presentations.  The game is played within each group, for the new careers, occupations, and job titles selected between Week 5 and Week 8, until each student has had at least two turns to present careers and choose cards.}

\handin{5}{8}

\vspace{4mm}

\noindent \hrulefill

\vspace{2mm}

\assignment{5}{8}{13}

\vspace{2mm}

\fordeadline{five weeks}
\begin{itemize}
\item  For \emph{\small{ALL}} of the remaining profiles that you have explored for the sake of our class for which you have not yet selected a proto-Guidelines card, do the following:
\begin{itemize}
\item  Select, from among the cards ``won'' during our in-class game \emph{and} the cards still in your hand, the ethical proto-guideline that you think is most relevant for someone in each selected career or position or job title.\ \ \onetoone
\item  For those career-card pairings, perform the Rank-1 tasks.
\end{itemize}
\item  For two of the career-card pairings for which you have performed the Rank-1 tasks but not the Rank-2 tasks, complete the Rank-2 tasks.
\item  For two of the career-card pairings for which you have performed the Rank-2 tasks but not the Rank-3 tasks, complete the Rank-3 tasks.
\item  For two of the career-card pairings for which you have performed the Rank-3 tasks, also complete the ``Rank-4'' tasks on the attached worksheet.
\item  Bring all worksheets and relevant proto-Guidelines cards to class in five weeks. \ \ \discuss
\end{itemize}
\emph{Summary:}  When you come to class five weeks from now, you will have
\begin{itemize}
\item  two career-card pairings analyzed all the way to Rank-4 tasks,
\item  two other career-card pairings analyzed all the way to Rank-3 tasks,
\item  two more career-card pairings analyzed with respect to Rank-2 tasks, and
\item  four more career-card pairings analyzed with respect to Rank-1 tasks.
\end{itemize}

\vspace{4mm}

\noindent \hrulefill

\vspace{4mm}

\begin{center}
\textbf{Rank-4 tasks (to be completed as homework)}
\end{center}

\vspace{4mm}

\noindent \emph{Career or position:}

\vspace{6mm}

\noindent Refer to your brainstorming and proposed actions from your Rank-1, Rank-2, and Rank-3 task work.  Respond to the following, giving each response in one to three sentences.

\vspace{6mm}

\noindent Propose a course of action that might be a reasonable option for the ethical mathematical 
practitioner in the situation.  (Proposing it does not imply that \emph{you} would choose it in the situation; it implies only that you are \emph{aware of it} as an option.)

\vspace{6mm}

\noindent Describe how the practitioner could defend the decision, specifically to those
with various opposing ethical frameworks or to those who incur costs from the proposed course 
of action.

\vspace{6mm}

\noindent Tell how and when the practitioner might reasonably reflect on the course of 
action---the ramifications of the action, implications for similar decisions in the 
future---and criteria that might be used.

\vspace{6mm}

\noindent \hrulefill

\vspace{4mm}

\noindent \textbf{{\small{SAMPLE}} class activity for Week 13.}  \emph{At the beginning of class during Week 13, students deliver their brief presentations and hand in their worksheets and cards.  They hand in, also, ``big picture'' worksheet.}

\emph{By the end of class, the instructor hands {\small{BACK}} each student's cards as collected over the course of the semester---including the ones handed in that day, in Week 13---and distributes the final assignment:}

\vspace{4mm}

\noindent \hrulefill

\vspace{2mm}

\assignment{6}{13}{of Finals}

\vspace{2mm}

\emph{Due at the Final Exam:}
\begin{enumerate}
\item  Carry out research in one of the following ways with the source described.
\begin{itemize}
\item  Read the list of Data Ethics Tenets, listed on page 4 and explored on pages 10--25, in the \emph{Federal Data Strategy Data Ethics Framework} at \texttt{https://resources.data.gov/assets/documents/} \linebreak \texttt{fds-data-ethics-framework.pdf}.  (These Tenets are required for all users of data who work for the Federal government of the United States.)  There are seven Data Ethics Tenets.

\item  Read \emph{The Toronto Declaration: Protecting the right to equality and non-discrimination in machine learning systems} prepared by Amnesty International and Access Now and available at \linebreak \texttt{https://www.accessnow.org/wp-content/uploads/2018/08/} \\ \texttt{The-Toronto-Declaration$\underline{\ }$ENG$\underline{\ }$08-2018.pdf}.  There are 49 articles (Articles 8 to 56) that are not part of the Introduction or the Conclusion.

\item  Find and read the professional ethics code or professional practice standards for some field or profession of your choice \emph{other than} mathematics, statistics, or computer science.

\end{itemize}
We will refer from now on to the ``chosen list of ethical standards'' and understand that we are referring to items from one of the works noted---Data Ethics Tenets, articles from the \emph{Toronto Declaration}, or points of another profession's ethics code or practice standards.

\item  Imagine that you hold a position with the Federal government that requires you to handle data according to the Data Ethics Tenets, or that you are working in artificial intelligence or a related field and wish to abide by the articles of the Toronto Declaration, or that you are a professional in another field whose code of ethics or practice standards you researched.

You have collected \emph{and played} a hand of ten proto-Guidelines cards.  Reflect on the items in the chosen list of ethical standards; select seven of them; and, for each of those seven standards, do the following.
\begin{enumerate}
\item  Identify the ethical mathematical practice proto-Guideline element, from among your ten cards, that, if observed consistently, will \emph{support} meeting that standard.  That is, for which of the proto-Guideline elements in your hand could you say, ``If I abide by these proto-Guideline elements, then they will help me to meet this other standard''?
\item  Determine whether there is any single ethical mathematical practice proto-Guideline element in your hand, or whether there is any subset of them, that, if observed consistently, would \emph{guarantee} meeting the standard.  That is, is there a subset of the proto-Guideline elements in your hand for which you could say, ``If I abide by these proto-Guideline elements, then I will absolutely have met this other standard''?
\end{enumerate}

\item  Write a short paper [\emph{specifications set by the instructor}] relating your responses to the previous questions and giving reasong for your answers.

\end{enumerate}

\vfill

\section*{APPENDIX B}

The pages that follow are printable versions of the proto-Guidelines cards.  We recommend printing the cards on white cardstock of density at least 176 grams per square centimeter.  Each sheet needs to be printed on twice: once for the card faces and once for the card backs.

\emph{The proto-Guidelines for Ethical Mathematical Practice were developed by Buell, Piercey, and Tractenberg~\cite{BPT2022} and are presented in Tractenberg, Piercey, and Buell~\cite{TPB2024PG} and~\cite{TPB2024}.  They are used here by permission of the authors.}

\pagebreak

\ \hspace{-8.5mm} \includegraphics[scale=0.90]{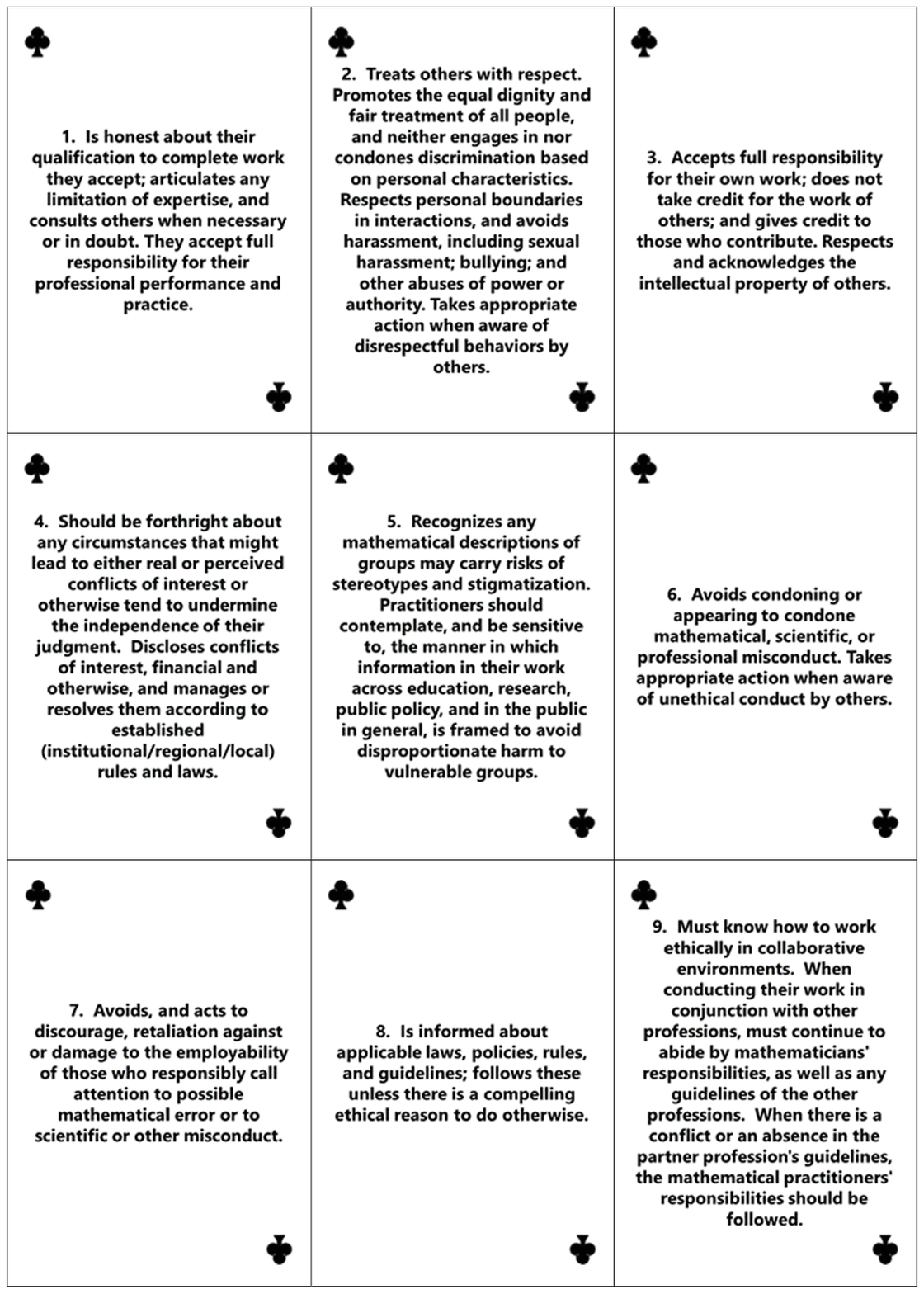}

\pagebreak

\ \hspace{-8.5mm} \includegraphics[scale=0.90]{EthicsDeckDraft220231221_2.eps}

\pagebreak

\ \hspace{-8.5mm} \includegraphics[scale=0.90]{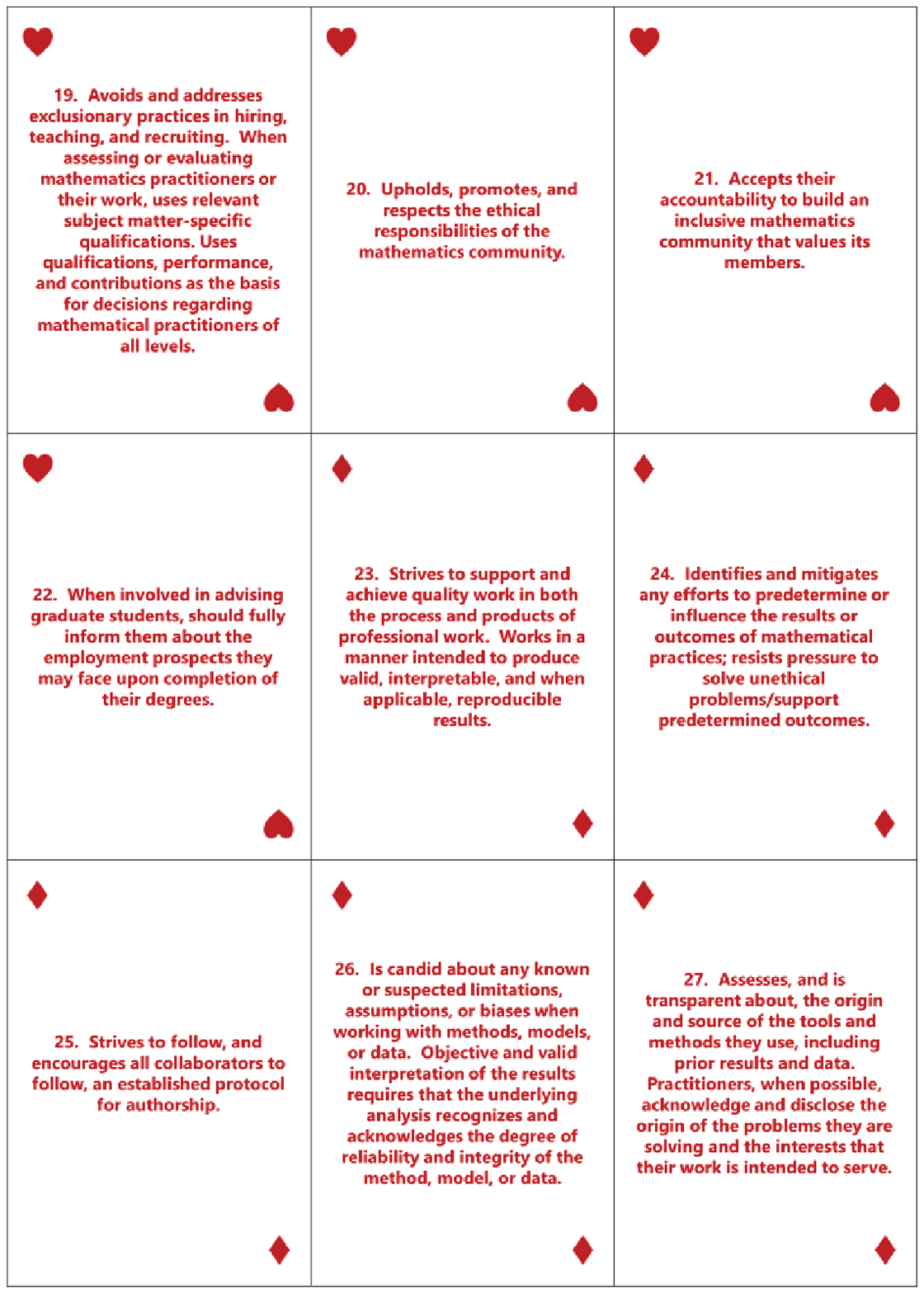}

\pagebreak

\ \hspace{-8.5mm} \includegraphics[scale=0.90]{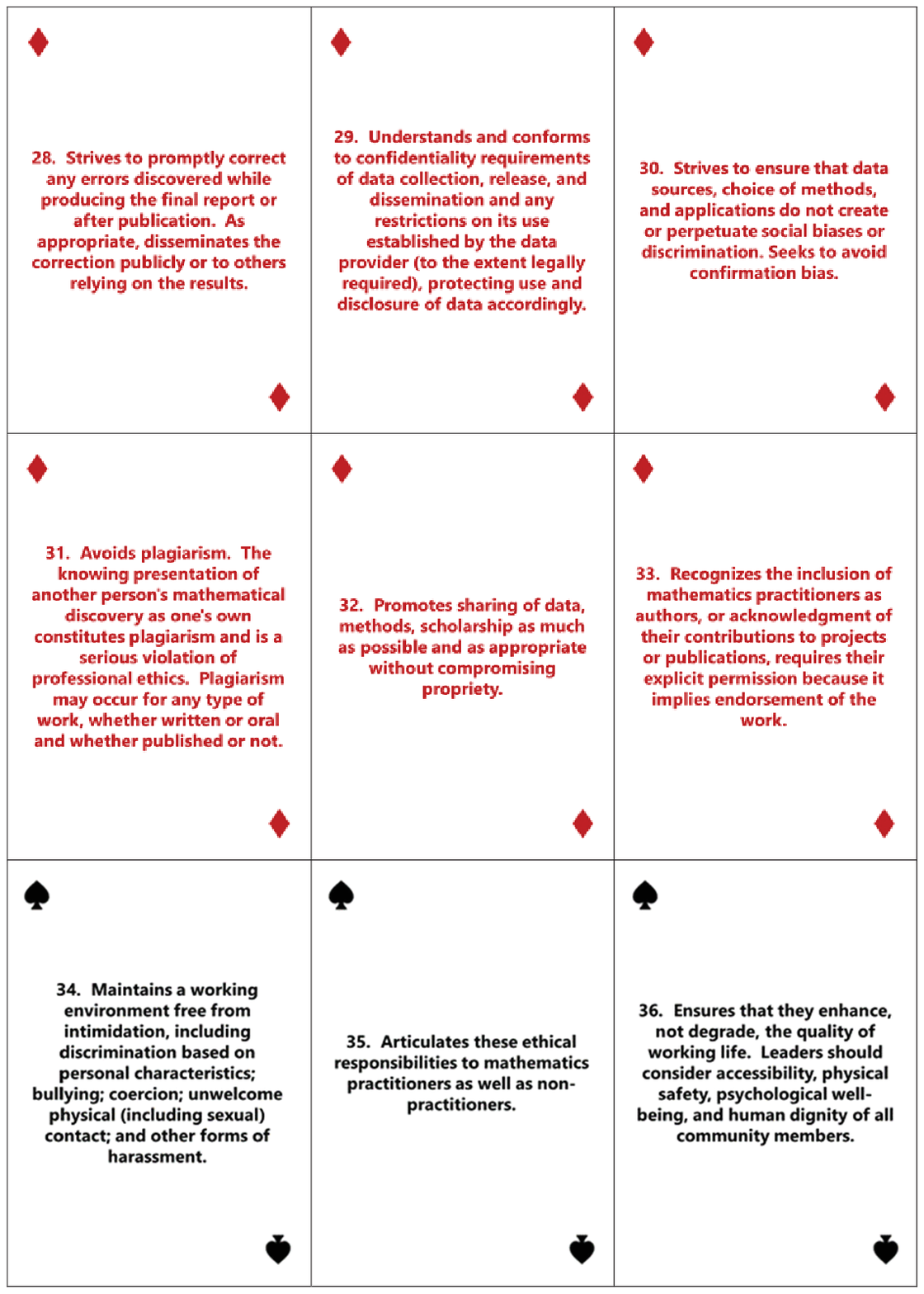}

\pagebreak

\ \hspace{-8.5mm} \includegraphics[scale=0.90]{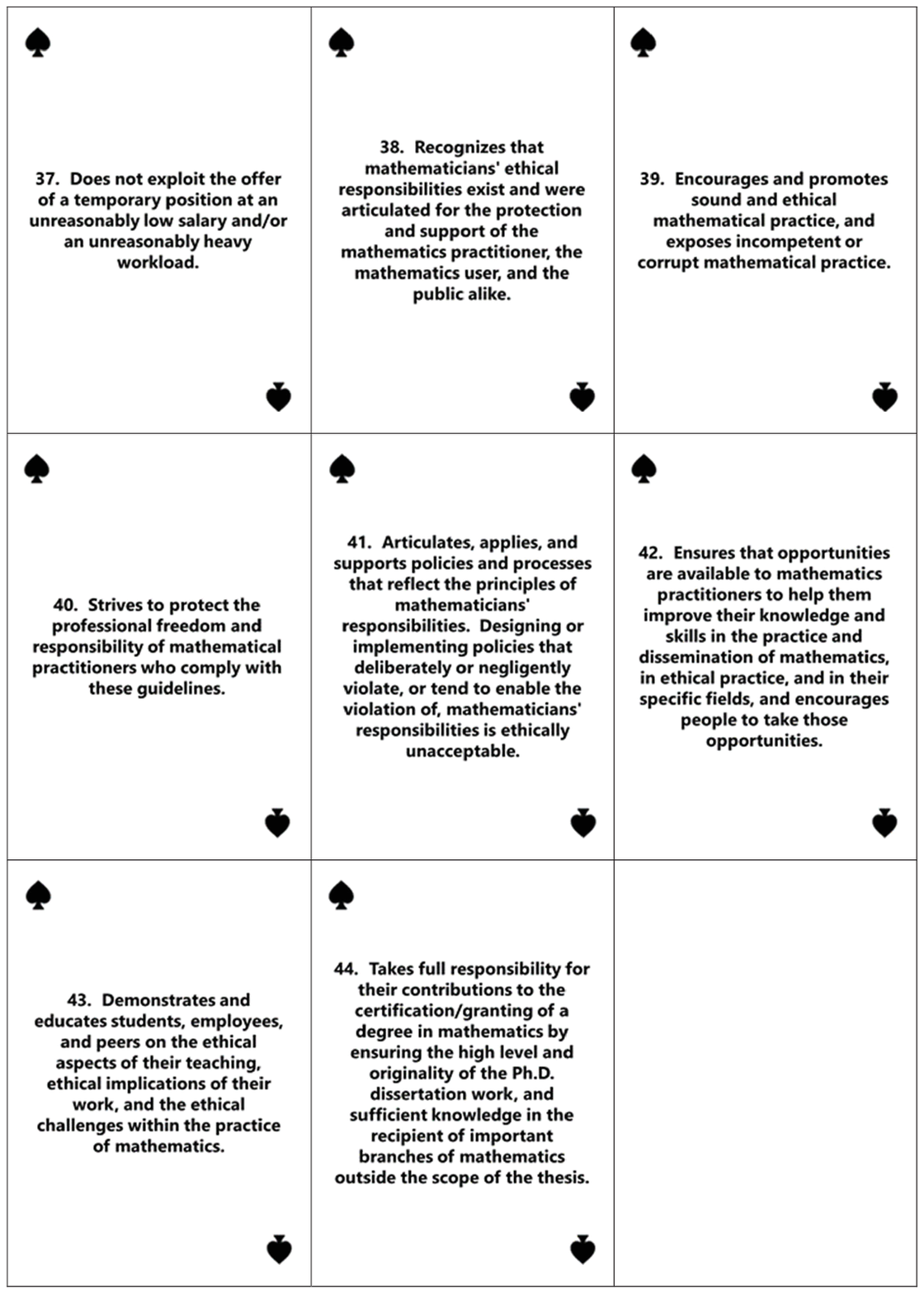}

\pagebreak

\ \hspace{-8.5mm} \includegraphics[scale=0.90]{EthicsDeckDraft220231221_6.eps}

\end{document}